\documentclass{article}  

\usepackage{graphicx}
\usepackage{amsthm}
\usepackage{amsfonts}
\usepackage{amsmath}
\usepackage{psfrag}

\newtheorem{theorem}{Theorem}[section]

\newtheorem*{propositionB}{Proposition 1.5}

\newtheorem{lemma}[theorem]{Lemma}
\newtheorem{proposition}[theorem]{Proposition}
\newtheorem{corollary}[theorem]{Corollary}
\theoremstyle{definition}
\newtheorem{definition}[theorem]{Definition}
\theoremstyle{remark}
\newtheorem{remark}[theorem]{Remark}
\newtheorem{problem}{Problem}

\numberwithin{equation}{section}
\numberwithin{figure}{section}

\newcommand{\hyperbolic}{\mathbb{H}}  
\newcommand{\complex}{\mathbb{C}} 
\newcommand{\integers}{\mathbb{Z}}   
\newcommand{\presentation}[2]{\langle \, #1  \, | \, #2 \,  \rangle}
\newcommand{\GammaH}{Cayley (G, X\cup \widetilde{\mathcal{H}}) }
\newcommand{\Hs}{\{ H_i \}_{i=1}^m}

\begin{document}

\title{Combination of Quasiconvex Subgroups of Relatively Hyperbolic Groups}
\author{Eduardo Mart\'inez-Pedroza \\ Department of Mathematics and Statistics\\ McMaster University \\ Hamilton, Ontario L8S 4K1.  Canada \\ emartinez@math.mcmaster.ca}
\date{}

\maketitle

\begin{abstract}
For relatively hyperbolic groups, we investigate conditions guaranteeing that the subgroup generated 
by two quasiconvex subgroups $Q_1$ and $Q_2$ is quasiconvex and isomorphic to $Q_1 \ast_{ Q_1 \cap Q_2} Q_2$.
Our results generalized known combination theorems for quasiconvex subgroups of word-hyperbolic groups. 
Some applications are presented.
\end{abstract}

\section{Introduction}\label{sec.introduction}

Relatively hyperbolic groups, originally introduced by M.\ Gromov ~\cite{Gr87}, have received a great deal of attention by group theorists after foundational works
by B.\ Farb ~\cite{Fa98} and B.H.\ Bowditch ~\cite{BO99} in the late nineties.  This class of groups includes many interesting subclasses. For instance,
limit groups which are an essential part of the theory of algebraic geometry over free groups ~\cite{Al05, DA03}, and geometrically finite Kleinian groups which contain fundamental groups of finite-volume hyperbolic manifolds.

If $G$ is a countable group and $\mathcal{H}$ is a collection of subgroups of $G$, the notion of relative hyperbolicity for the pair $(G, \mathcal{H})$  has been defined by 
different authors ~\cite{BO99, DS05, Fa98, Gr87, GM06, HK08, Os06, Ya99}.  All these definitions are equivalent when the group $G$ and the subgroups in $\mathcal{H}$
are finitely generated ~\cite{ DS05, Gr87, Os06, Sz98, Ya99}. A precise definition is provided in the next section. 
When a pair  $(G, \mathcal{H})$ satisfies the relative hyperbolicity condition we say that the group $G$ is hyperbolic relative to $\mathcal{H}$, and 
when the collection $\mathcal{H}$ is fixed we just say that the group $G$ is relatively hyperbolic. 

For a group $G$ hyperbolic relative to a collection of subgroups $\mathcal{H}$, 
the quasiconvex subgroups are the natural subgroups to study when considering a relatively hyperbolic group as a geometric object. 
Different notions of ~\emph{relative quasiconvexity} for subgroups of $G$ were introduced by D.\ Osin and F.\ Dahmani ~\cite{ DA03, Os06} and 
recently C.\ Hruska has proved the equivalence of these definitions ~\cite{HK08}.

We are interested in the following problem.
\begin{problem}\label{main.problem}
Let $G$ be a relatively hyperbolic group, and suppose that $Q$ and $R$ are quasiconvex subgroups of $G$. 
Consider the natural homomorphism
\begin{eqnarray*}  \rho: Q \ast_{Q\cap R} R \longrightarrow G  ,\end{eqnarray*}
which has image the subgroup $\langle Q\cup R\rangle$.
\begin{enumerate}
\item (Algebraic Structure.)   When is $\rho$ injective?
\item (Geometric Structure.)  When is the image of $\rho$ a quasiconvex subgroup?
\end{enumerate}
\end{problem}

\subsection{Main Results}
Let $G$ be a group generated by a finite set $X$ and hyperbolic relative to a collection of subgroups $\mathcal{H}$. 
A subgroup of $G$ is called \emph{parabolic} if can be conjugated into one the subgroups in $\mathcal{H}$. For an element $g\in G$, $|g|_X$ denotes its distance
from the identity element in the word metric induced by $X$ on $G$

\begin{theorem}[Quasiconvex-Parabolic amalgamation]
 \label{thm:main_intro}
For any relatively quasiconvex subgroup $Q$ and any maximal parabolic subgroup $P$ of $G$,  there is constant $C=C(Q,P) \geq 0$ with the following property.
If $R$ is a subgroup of $P$ such that
\begin{enumerate}
\item $Q\cap P \leq R$, and
\item $|g|_X \geq C$ for any $g \in  R \setminus  Q$,
\end{enumerate}
then the natural homomorphism
\begin{eqnarray*}  Q \ast_{Q\cap R} R \longrightarrow G   \end{eqnarray*}
is injective with image a relatively quasiconvex subgroup. \newline
Moreover, every parabolic subgroup of $\langle Q\cup R \rangle < G$ is either conjugate to a subgroup of $Q$ or a subgroup of $R$ in $\langle Q\cup R \rangle$.
\end{theorem}

\begin{theorem}[Quasiconvex-Quasiconvex amalgamation]  \label{thm:main_intro_2}
For any pair of  relatively quasiconvex subgroups $Q_1$ and $Q_2$, and any maximal parabolic subgroup $P$ such that $R=Q_1\cap P = Q_2 \cap P$, 
there is a constant $C=C(Q_1, Q_2, P) \geq 0$ with the following property.\newline
If $h \in P$ is such that
\begin{enumerate}
\item $hRh^{-1} = R$, and
\item $|g|_X \geq C$ for any $g \in  RhR $,
\end{enumerate}
then the natural homomorphism
\begin{eqnarray*}  Q_1 \ast_R hQ_2h^{-1} \longrightarrow G   \end{eqnarray*}
is injective and its image is a relatively quasiconvex subgroup.\newline
Moreover, every parabolic subgroup of $\langle Q_1 \cup hQ_2h^{-1} \rangle < G$ is either conjugate to a subgroup of $Q_1$ or $hQ_2h^{-1}$ in $\langle Q_1 \cup hQ_2h^{-1} \rangle$.
\end{theorem}

\subsection{History and Motivation}
This work is motivated by known combination theorems for quasiconvex subgroups of word-hyperbolic groups.
In ~\cite{Gr87}, M.\ Gromov stated that in a torsion free word-hyperbolic group any infinite index quasiconvex subgroup is a free factor of a larger quasiconvex subgroup.
Gromov's ideas were developed by G.N.\ Arzhantseva in ~\cite{Az01}. More general combination theorems for quasiconvex subgroups of word-hyperbolic groups 
were stated and proved by R.\ Gitik  in ~\cite{Gi99}. For relatively hyperbolic groups with discrete representations in $Isom(\hyperbolic^n)$, 
recent results by M.\ Baker and D.\ Cooper correspond to combination of quasiconvex subgroups ~\cite{BC05}. 

The Klein-Maskit Combination Theorems for Kleinian groups ~\cite{Maskit} are another motivation for our work.
In particular, the following example whose details can be found in ~\cite{LR06}: if $G_1$ and $G_2$ are two lattices of $PSL(2, \complex)$
and $R$ is a maximal parabolic subgroup of both, then for a ``sufficiently complicated" parabolic $h$ centralizing $R$,  
the natural homomorphism from $Q_1 \ast_{R} h Q_2h^{-1}$ into $PSL(2, \complex)$ is injective. 
This technique has been used by D.\ Cooper, D.\ Long, and A.\ Reid to double quasiFuchsian subgroups along parabolic subgroups in hyperbolic manifold groups, 
producing essential closed surfaces in cusped hyperbolic manifolds  ~\cite{CL01, CLR97}.
Corollary~\ref{cor:double_quasiconvex} illustrates this technique in the context  of relatively hyperbolic groups. 

Another motivating result, with a similar statement to Theorem~\ref{thm:main_intro_2}, is 
a combination theorem for Veech subgroups of the mapping class group by  C.\ Leininger and A.\ Reid in ~\cite{LR06}. 
This result was used to construct subgroups of the mapping class group isomorphic to the fundamental group of a closed surface.
The mapping class group is weakly relatively hyperbolic  ~\cite{MM99}, but is not strongly relatively hyperbolic ~\cite{AAS05, BDM05}. 

The main motivation of this work is an outstanding question by M.\ Gromov of whether every one-ended word-hyperbolic group contains a subgroup isomorphic to the fundamental group of a closed surface.  
In ~\cite{CL01}, D.Cooper and D.Long produce surface subgroups in word-hyperbolic groups arising as the fundamental groups of Dehn fillings  of finite volume hyperbolic manifolds, starting with the construction of 
quasiconvex subgroups with particular structures in the finite volume hyperbolic manifold group. We aim to explore Gromov's question on particular classes of word-hyperbolic groups which arise as algebraic Dehn fillings of relatively hyperbolic groups. 
The notion of algebraic Dehn filling has been studied by D.\ Groves and J.\ Manning, and independently by D.\ Osin ~\cite{GM06, Os06-1}. 
The main results of this paper are part of this program.

\subsection{Other results on Quasiconvex Subgroups}
Let $G$ be a hyperbolic group relative to a collection of subgroups $\mathcal{H}$ with finite generating set $X$.

\begin{proposition}\label{prop:IntersectionQuasiconvexGroups}
Let $Q$ and $R$ be relatively quasiconvex subgroups of $G$. Then $Q\cap R$ is a relatively quasiconvex subgroup of $G$.
\end{proposition}
\begin{remark}
In the case of word-hyperbolic groups, Proposition ~\ref{prop:IntersectionQuasiconvexGroups} was originally proved by H.\ Short. 
C.\ Hruska has independently proved this property without assuming that the ambient group is finitely generated ~\cite{HK08}.
In this generality, D.\ Osin stated the same result in ~\cite{Os06}. 
\end{remark}

\begin{proposition}\label{prop:ParabolicClasses}
Let $Q$ be a $\sigma$-quasiconvex subgroup of $G$. The number of infinite maximal parabolic subgroups of $Q$ up to conjugacy in $Q$ is finite.
\end{proposition}
\begin{remark}
Proposition ~\ref{prop:ParabolicClasses} appears as ~\cite[Theorem 9.1]{HK08} and it is proved using the dynamical characterization of quasiconvexity. 
Here we present a conceptually simpler proof using arguments on Cayley graphs.
\end{remark}

\subsection{Sample Applications}

Special attention has been given to relatively hyperbolic groups with peripheral structure consisting of abelian or virtually abelian subgroups.
In this setting Theorems~\ref{thm:main_intro} and~\ref{thm:main_intro_2} can be used to construct quasiconvex subgroups with particular structures.

Let $G$ be hyperbolic relative to a collection of free abelian subgroups $\mathcal{H}$.

Given a subgroup $R$ of a group $Q$, the amalgamated free product of $k$ copies of $Q$ along $R$ is denoted by $\Delta_m (Q, R)$. 
When $m=2$, $\Delta_2 (Q, R)$ is called {\it the double of $Q$ along $K$}. Doubling a group along a subgroup has been used in different contexts in 
group theory and geometric topology; for example to produce groups with interesting finiteness properties ~\cite{Sta63}, or to produce surface subgroups in finite volume hyperbolic 3-manifold groups ~\cite{CLR97}.

\begin{corollary}[Doubling Quasiconvex along Parabolics]\label{cor:double_quasiconvex}
Let $Q$ be a relatively quasiconvex subgroup and let $P$ be a maximal parabolic subgroup of $G$.
If 
\[rank_{\integers} (Q\cap P) < rank_{\integers} (P),\]
then there exists a quasiconvex subgroup isomorphic to $\Delta_k (Q, Q\cap H)$ for any positive integer $k$.
\end{corollary}

A quasiconvex subgroup $R$ of $G$ is called {\it fully quasiconvex} if for any parabolic subgroup $P < G$, the subgroup $Q\cap P$ is finite or of finite index in $P$. 
Fully quasiconvex subgroups have appeared in the work of F.\ Dahmani ~\cite{DA03} where is shown, under some hypothesis, that the 
combination of relatively hyperbolic groups along fully quasiconvex subgroups is a relatively hyperbolic group. 

\begin{corollary}[Fully quasiconvex amalgamams] \label{cor:fully_quasiconvex}  
Let $Q$ be a relatively quasiconvex subgroup.  Then there exists a fully quasiconvex subgroup $R$ which splits over $Q$. 
\end{corollary}

\subsection{Outline of the Paper}
The paper is organized as follows. Section~\ref{sec.hyperbolicity} introduces background and notation. 
Section~\ref{sec.quasigeodesics} states and proves a proposition about quasi-geodesics which complements a result by C.\ Dru\c tu and M.\ Sapir in \cite{DS05}. 
Section~\ref{sec.quasiconvexity} recalls the notion of relatively quasiconvex subgroup and the proofs of Propositions ~\ref{prop:IntersectionQuasiconvexGroups} and~\ref{prop:ParabolicClasses} are explained.
Section~\ref{sec.mainproofs} consists of the proofs of the main results and the applications.  The last section of the paper consists of future research directions.

\subsection{Acknowledgements} 
The results of this project are part of the author's PhD dissertation at the University of Oklahoma.
The author thanks his academic advisor Noel Brady for his guidance and encouragement, and Boris Apanasov, Ara Basmajian, Max Forester, Christopher Leininger, Pallavi Dani, Darryl McCullough, Krishnan Shankar, and Stephen Weldon for helpful comments during this work. 
The author also thanks the referee for several useful comments.
This project was partially supported by NSF grant no. DMS-0505707 and the Department of Mathematics at the University of Oklahoma through a Foundation Fellowship. 

\section{Relatively Hyperbolic Groups} ~\label{sec.hyperbolicity}

The aim of this section is to introduce notation and to define relatively hyperbolicity for finitely generated groups.
The definition presented below is equivalent to the one given by D.\ Osin in ~\cite{Os06}; the equivalence follows directly from  ~\cite[Theorems 3.23 and 6.10]{Os06}.

\subsection{Preliminaries}
We follow closely the notation and conventions of the paper ~\cite{Os06}.
Let $G$ be a group and $A \subset G$ a generating set closed under inverses.
The {\it Cayley Graph of the group $G$ with
respect to $A$}, which is denoted by $Cayley (G, A)$, is the oriented graph with vertex set $G$ and edge set $G \times A$, where
an edge $e=(g,a)$ goes from the vertex $g$ to the vertex $ga$ and has label $Label(e) = a$.
Let $p = e_1e_2 \dots e_k$ be a combinatorial path in $Cayley(G,A)$.
The initial and the terminal vertices of $p$ are denoted by
$p_-$ and $p_+$ respectively, the label $Label(p)$ of $p$ is the word
$Label (e_1) Label (e_2) \dots Label (e_k)$ in the alphabet $A$, and
the length $l(p)$ of $p$ is the number of edges in $p$.
The concatenation of the combinatorial paths $p$ and $q$ such that $p_+=q_-$ is denoted by $p q$.
The (word) length $|g|_A$ of an element $g \in G$ is the length of a shortest
combinatorial path in $Cayley(G, A)$ from $1$ to $g$. This defines a left invariant
metric on the group $G$ defined by $dist_A (f, g) = |f^{-1}g|_A$.

A geodesic metric space $(\Gamma, d)$ is a {\it $\delta$-Gromov hyperbolic} for
some $\delta \geq 0$ if for any geodesic triangle $\Delta$, every side is contained in the
open $\delta$-neighborhood of the other two sides. A rectifiable path $p$ in $(\Gamma , d)$ is a
{\it $(\lambda, c)$-quasi-geodesic } for some $\lambda \geq 1$ and $c\geq 0$ if for any subpath $q$ of $p$
\[ l(q) \leq \lambda d(q_-, q_+) + c.\]

\subsection{Definition of Relative hyperbolicity.}
Let $G$ be a group, $\mathcal{H}=\Hs$ be a collection of subgroups of $G$, and $X$ be a symmetric finite generating set for $G$. 
Denote by $\widetilde{\mathcal{H}}$ the disjoint union
\begin{eqnarray*}  \widetilde{\mathcal{H}} = \bigsqcup_{i=1}^m (\widetilde{H_i} \setminus \{1\})   \end{eqnarray*}
where $\widetilde{H_i}$ is a copy of $H_i$.

\begin{definition}[Weak Relative Hyperbolicity]
The pair $(G, \Hs)$ satisfies the \emph{weakly relative hyperbolicity} condition if there is an integer $\delta \geq 0$ such that the Cayley graph $\GammaH$ is a $\delta$-hyperbolic metric space.
\end{definition}

\begin{definition}[ D.\ Osin ~\cite{Os06}.  $\mathcal{H}$-components, connected and isolated, backtracking, phase vertices, and $k$-similar paths ]

Let $q$ be a combinatorial path in the Cayley graph $\GammaH$. Subpaths of $q$ with at least one edge are called non-trivial.
An {\it $H_i$-component} of $q$ is a maximal non-trivial subpath $s$ of $q$ with $Label(s)$ a word in the alphabet $\widetilde{H_i} \setminus \{1\}$. 
When we don't need to specify the index $i$, we will refer  $H_i$-components as $\mathcal{H}$-components.

Two $\mathcal{H}$-components $s_1$, $s_2$ of $q$ are {\it connected} if the vertices of $s_1$ and $s_2$ belong to the same left coset of $H_i$ for some $i$.
Equivalently, the $\mathcal{H}$-components $s_1$ and $s_2$ are connected if $Label(s_1)$ and $Label(s_2)$ are words in the alphabet $\widetilde{H_i}$ for some $i$,
there exists a path $c$ which connects a vertex of $s_1$ and a vertex of $s_2$, and $Label(c)$ is a word in the alphabet $\widetilde{H_i}$.
An $\mathcal{H}$-component $s$ of $q$ is {\it isolated} if it is not connected to a different $\mathcal{H}$-component of $q$.

The path $q$ is {\it without backtracking} if every $\mathcal{H}$-component of $q$ is isolated. 
A vertex $v$ of $q$ is called {\it phase} if it is not an inner vertex of an $\mathcal{H}$-component $s$ of $q$. Two paths $p$ and $q$ in $\GammaH$ are {\it $k$-similar} if
\[ max\{d_{X}(p_-,q_-), d_{X}(p_+, q_+)  \} \leq k.\]
\end{definition}
\begin{remark}
Every geodesic path in $\GammaH$ is  without backtracking and all its vertices are phase.
\end{remark}

\begin{definition}[Bounded Coset Penetration (BCP)] \label{BCP}
The pair $(G, \Hs)$ satisfies the BCP property if for any $\lambda \geq 1$, $c \geq 0$, $k\geq 0$, there exists an integer $\epsilon(\lambda, c, k) > 0$ such that for any two $k$-similar $(\lambda, c)$-quasi-geodesics in $\GammaH$ without backtracking $p$ and $q$ the following holds.
\begin{enumerate}
\item[(i.)] The sets of phase vertices of $p$ and $q$ are contained in the closed
$\epsilon(\lambda, c, k)$-neighborhoods (with respect to the metric $dist_X$) of each other.
\item[(ii.)] Suppose $s$ is an $\mathcal{H}$-component of $p$ such that $dist_X(s_-, s_+) > \epsilon(\lambda, c, k)$; then
there exists an $\mathcal{H}$-component $t$ of $q$ which is connected to $s$.
\item[(iii.)] Suppose $s$ and $t$ are connected $\mathcal{H}$-components of $p$ and $q$ respectively.
Then $max\{ dist_X(s_-, t_-), dist_X(s_+, t_+) \} \leq  \epsilon(\lambda, c, k)$.
\end{enumerate}
\end{definition}
\begin{remark}
Our definition of the BCP property corresponds to the conclusion of Theorem 3.23 in \cite{Os06}.
\end{remark}
\begin{definition}[Relative Hyperbolicity] \label{defn:rel_hyp}
The pair $(G, \Hs)$ satisfies the \emph{relative hyperbolicity} condition if the group $G$ is weakly hyperbolic relative to $\Hs$ and the pair $(G,\Hs)$ satisfies the Bounded Coset Penetration property.
If $(G, \Hs)$ satisfies the relative hyperbolicity condition then we said that group $G$ is hyperbolic relative to $\Hs$; if there is no ambiguity, we just said that the group $G$ is relatively hyperbolic.
\end{definition}
\begin{remark}
The stated definition of relative hyperbolicity is equivalent to  \cite[Definition 2.35]{Os06} for finitely generated groups; 
the equivalence follows directly from \cite[Theorems 3.23 and 6.10]{Os06}.
\end{remark}

For the rest of this section, 
let $G$ be a group hyperbolic relative to a collection of subgroups $\mathcal{H}$ and let $X$ be a symmetric finite generating set of $G$. 
\begin{figure}[ht]
\begin{center}
\psfragscanon{
\psfrag{$p$}{$p$}
\psfrag{$q$}{$q$}
\psfrag{$p_-$}{$p_-$}
\psfrag{$p_+$}{$p_+$}
\psfrag{$q_-$}{$q_-$}
\psfrag{$q_+$}{$q_+$}
\psfrag{$g_1H_i$}{$g_1H_i$}
\psfrag{$g_2H_j$}{$g_2H_j$}
\includegraphics[width=0.8\textwidth]{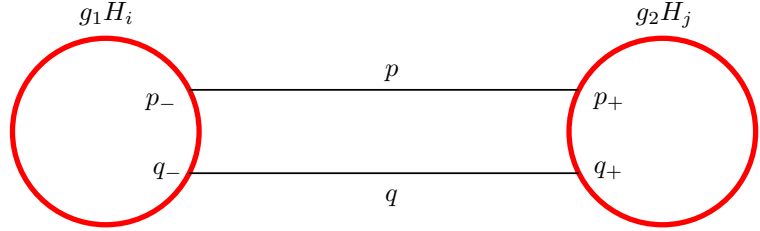}
}
\end{center}
\caption{ Two geodesics connecting different left cosets $g_1H_i$ and $g_2H_j$ are $\epsilon(1,4,0)$-similar, and their lengths differ by at most two (Corollary~\ref{cor:TubeLemma}).}
\end{figure}
The following corollary is a direct consequence of Definition ~\ref{defn:rel_hyp}.
\begin{corollary}\label{cor:TubeLemma}
Let $g_1H_i$ and $g_2H_j$ be different left cosets.
For any pair of geodesics $p$ and $q$ in $\GammaH$ such that $p_-, q_- \in g_1H_i$, $p_+, q_+ \in g_2H_j$, and neither $p$ nor $q$ have more than one vertex in $g_1H_i$ or $g_2H_j$, the following holds.
\begin{enumerate}
\item $l(q) \leq l(p) + 2$, and
\item $q$ and $p$ are $\epsilon(1,4,0)$-similar.
\end{enumerate}
\end{corollary}
\begin{proof}
Consider the path $r = c_1  p  c_2$ in $\GammaH$ where $c_1$ is an edge connecting
$q_-$ and $p_-$, and $c_2$ is an edge connecting $p_+$ and $q_+$. Notice that $r$ is a $(1,4)$-quasi-geodesic
in $\GammaH$ and that $q$ and $r$ are $0$-similar.
The BCP-property implies
\begin{eqnarray*} max \{ dist_X( (c_1)_-, (c_1)_+ ),  dist_X( (c_2)_-, (c_2)_+ )    \} \leq \epsilon(1,4,0)  .  \end{eqnarray*}
\end{proof}

\section{Quasi-geodesics}~\label{sec.quasigeodesics}

Let $G$ be a group generated by a finite set $X$, $\Hs$ a collection of subgroups of $G$, and suppose that
$G$ is hyperbolic relative to  $\Hs$.  Any geodesic $p$ in $\GammaH$ can be decomposed as
\[ p = r_1  s_1  \dots  r_k  s_k \]
where each $r_i$ is a geodesic, and each $s_i$ is an isolated $\mathcal{H}$-component.

In this section, we investigate paths with the above type of decomposition and estimate 
quasi-geodesic constants. The main result of the section is the following.
\begin{figure}[ht]
\begin{center}
\psfragscanon{
\psfrag{$r_1$}{$r_1$}
\psfrag{$s_1$}{$s_1$}
\psfrag{$r_2$}{$r_2$}
\psfrag{$s_2$}{$s_2$}
\psfrag{$r_3$}{$r_3$}
\psfrag{$s_3$}{$s_3$}
\psfrag{$r_4$}{$r_4$}
\psfrag{$s_4$}{$s_4$}
\includegraphics[width=0.95\textwidth]{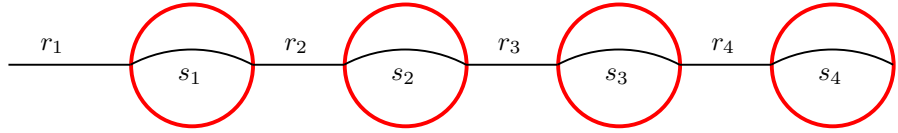}
}
\end{center}
\caption{Polygonal path $p$ in Proposition~\ref{prop:nicepath}.}
\end{figure}
\begin{proposition}\label{prop:nicepath}
There are constants $\eta$ and $\lambda_0$ with the following properties.
If $p$ is a path in $\GammaH$  such that
\begin{enumerate}
\item $p = r_1  s_1  \dots  r_k  s_k$ where each $r_i$ and $s_i$
are geodesic paths in $\GammaH$,
\item the subpath $s_i$ of $p$ is an $\mathcal{H}$-component of $p$ for each $i$,
\item $dist_X((s_i)_- ,(s_i)_+) \geq  \eta$ for each $i$, and
\item the $\mathcal{H}$-components $s_i$ and $s_{i+1}$ of $p$ are not connected for each $i$.
\end{enumerate}
Then $p$ is a $(\lambda_0,0)$-quasi-geodesic without backtracking and with different endpoints.
\end{proposition}

\begin{remark}
The strength of this result is that the constant $\lambda_0$ is independent of the number 
of segments of the path.  This complements a similar result by C.\ Dru\c tu and M.\ Sapir ~\cite[Lemma 8.12]{DS05}.
Assuming that the lengths of the segments $r_i$ are larger than a fixed constant $l$,  
their result estimates quasi-geodesic constants depending on $l$ and the number of segments $2k$.
\end{remark}

The rest of the section consists of two parts. First, a result by D.\ Osin about polygons in $\GammaH$  
is recalled and a corollary is stated. The second part consists of a series of lemmas
and the proof of Proposition~\ref{prop:nicepath}.

\subsection{Osin's Result about Polygons}
The following proposition is a stronger version of
the Bounded Coset Penetration property. It is a central part of D.\ Osin's work in ~\cite{Os06-1}.
\begin{proposition}[D. Osin ~\cite{Os06-1} ] \label{thm:n-gon}
There exists a constant $D>0$ satisfying the following condition.
If $\mathcal{P} = p_1  p_2  \dots  p_n$ is an
$n$-gon in $\GammaH$ and $S \subset \{p_1, \dots , p_n\}$ such
that:
\begin{enumerate}
\item each side $p_i \in S$ is an isolated $\mathcal{H}$-component of $\mathcal{P}$, and
\item each side $p_i \not \in S$ is a geodesic path.
\end{enumerate}
Then
\[   \sum_{p \in S} dist_X( p_- , p_+ )  \leq  Dn \]
\end{proposition}
The next Corollary is used in the proof of the main result of the section. 
\begin{figure}[ht]
\begin{minipage}[c]{.3\linewidth}
\begin{center}
\psfragscanon{
\psfrag{$s_-$}{$s_-$}
\psfrag{$s_+$}{$s_+$}
\includegraphics[width=0.7\textwidth]{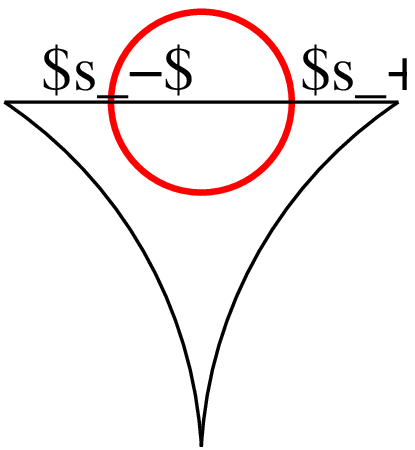}
}
\end{center}
\end{minipage}
\begin{minipage}[c]{.3\linewidth}
\begin{center}
\psfragscanon{
\psfrag{$s_-$}{$s_-$}
\psfrag{$s_+$}{$s_+$}
\psfrag{$t_-$}{$t_-$}
\psfrag{$t_+$}{$t_+$}
\includegraphics[width=0.7\textwidth]{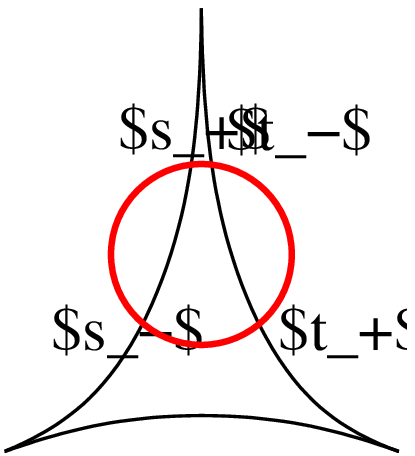}
}
\end{center}
\end{minipage}
\begin{minipage}[c]{.3\linewidth}
\begin{center}
\psfragscanon{
\psfrag{$s_-$}{$s_-$}
\psfrag{$s_+$}{$s_+$}
\psfrag{$t_-$}{$t_-$}
\psfrag{$t_+$}{$t_+$}
\psfrag{$u_-$}{$u_-$}
\psfrag{$u_+$}{$u_+$}
\includegraphics[width=0.7\textwidth]{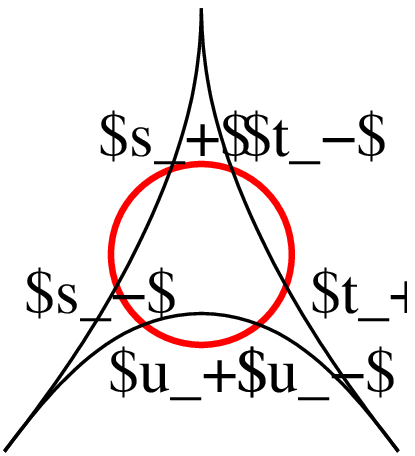}
}
\end{center}
\end{minipage}
\caption{The three cases of Corollary~\ref{cor:TraingleBigThm}.}
\end{figure}
\begin{corollary}\label{cor:TraingleBigThm}
There exists a constant $\tau \geq 0$ with the following property.
Let $\Delta = p  q  r$ be a triangle whose sides $p$, $q$, $r$ are geodesics
in $\GammaH$.
If $s$ is a $\mathcal{H}$-component of $p$, then either
\begin{enumerate}
\item $s$ is an isolated $\mathcal{H}$-component of the cycle $\Delta$ and
$dist_X(s_-,s_+) \leq \tau$;
\item $s$ is connected to only one $\mathcal{H}$-component $t$ of the concatenated path $q r$ and
 \[ max\{ dist_X(s_- ,t_+ ), dist_X(s_+ ,t_- ) \} \leq \tau ;\]
\item $s$ is connected to $\mathcal{H}$-components $t$ and $u$ of $q$ and $r$ respectively,  and
\[ max\{ dist_X(s_+,t_-), dist_X(t_+,u_-), dist_X(u_+,s_-) \} \leq \tau .\]
\end{enumerate}
\end{corollary}
\begin{proof}
In the first case, consider $\Delta$ as a $5$-gon with $s$ as one of its sides. Then
Proposition ~\ref{thm:n-gon} implies that $ dist_X(s_-, s_+) \leq 5D .$

In the second case, we consider the case $t$ is an
$\mathcal{H}$-component of $q$. Decompose the paths $p$ and $q$ as
$p = p_1  s  p_2$ and $q = q_1  t  q_2$.
Since $s$ and $t$ are connected $\mathcal{H}$-components, there
are edges $c_1$ connecting $t_-$ and $s_+$, and $c_2$ connecting
$s_-$ and $t_+$. Considering the $3$-gon $p_2  q_1 
c_1$ and the $4$-gon $p_1  c_2  q_2  r$, Proposition
~\ref{thm:n-gon} implies
\[ max\{ dist_X(s_- ,t_+ ), dist_X(s_+ ,t_- ) \} \leq 4D .\]

For the third case, an analogous argument shows that
\[ max\{ dist_X(s_+,t_-), dist_X(t_+,u_-), dist_X(u_+,s_-) \} \leq 3D .\]
To finish the proof define $\tau = 5D.$
\end{proof}
\begin{remark}
An equivalent result to Corollary~\ref{cor:TraingleBigThm} appears as ~\cite[Proposition 8.16]{DS05}.
\end{remark}

\subsection{Proof of Proposition~\ref{prop:nicepath}}

\begin{figure}[ht]
\begin{center}
\psfragscanon{
\psfrag{$p$}{$p$}
\psfrag{$q$}{$q$}
\psfrag{$r_1$}{$r_1$}
\psfrag{$s_1$}{$s_1$}
\psfrag{$t_1$}{$t_1$}
\psfrag{$r_2$}{$r_2$}
\psfrag{$s_2$}{$s_2$}
\psfrag{$t_2$}{$t_2$}
\psfrag{$r_3$}{$r_3$}
\psfrag{$s_3$}{$s_3$}
\psfrag{$t_3$}{$t_3$}
\psfrag{$r_4$}{$r_4$}
\psfrag{$s_4$}{$s_4$}
\psfrag{$t_4$}{$t_4$}
\includegraphics[width=0.95\textwidth]{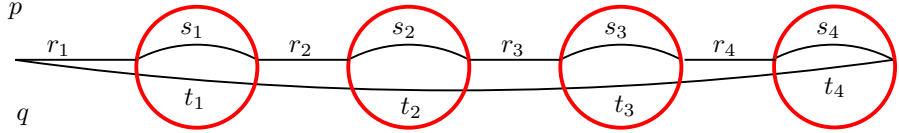}
}
\end{center}
\caption[Lemma~\ref{lemma:mainproposition}]
{The path $p$ and a geodesic $q$ connecting the endpoints of $p$.}
\label{fig:LemmaNicePath}
\end{figure}

We will see that $\lambda_0 = 3$ works, and will define a lower bound for $\eta$ during the
course of the proof. The argument consists of three lemmas.
\begin{lemma} \label{lemma:mainproposition}
Let $p$ be a path satisfying the hypothesis of Proposition ~\ref{prop:nicepath}, and let 
$q$ be a geodesic in $\GammaH$ connecting the endpoints of $p$.
For each $i$, $q$ contains an  $\mathcal{H}$-component $t_i$ satisfying :
\begin{enumerate}
\item $s_i$ and $t_i$ are connected $\mathcal{H}$-components, and
\item $s_i$ and $t_i$ are $\epsilon(1,4,0)$-similar.
\end{enumerate}
In particular, the endpoints of $p$ are different. (See Figure~\ref{fig:LemmaNicePath}.)
\end{lemma}
\begin{proof}
We argue by induction on $k$. 
Suppose $k=1$. Since $l(s_1)=1$ it follows that $r_1  s_1$ is a $(1,2)$-quasi-geodesic in $\GammaH$.
By our choice of $\eta$,
\begin{equation} dist_X((s_1)_-,(s_1)_+) \geq \eta > \epsilon (1,2,0).\end{equation}
The BCP-property implies that any geodesic  in $\GammaH$ connecting the endpoints of
$r_1  s_1$ has a $\mathcal{H}$-component which is connected to $s_1$ and is
$\epsilon(1,2,0)$-similar to $s_1$.

Suppose $k>1$. Consider the subpaths of $p$:
\begin{eqnarray*} p_1 = r_1  s_1  \dots  r_{k-1}  s_{k-1},  &  \  \text{and}  \  & p_2 = r_k  s_k .  \end{eqnarray*}
Let $q$, $q_1$, and $q_2$ be geodesics in $\GammaH$ connecting the endpoints of $p$, $p_1$ and $p_2$ respectively.

{\it Claim 1.} $q_1$ has no $\mathcal{H}$-component connected to the $\mathcal{H}$-component $s_k$. \newline
By induction hypothesis $q_1$ has a $\mathcal{H}$-component $u$ which is connected and $\epsilon(1,4,0)$-similar to $s_{k-1}$. 
By the triangle inequality and our choice of $\eta$, 
\begin{eqnarray}  
			 dist_X(u_-, u_+) & \geq & dist_X((s_{k-1})_-,(s_{k-1})_+) - 2 \epsilon(1,4,0)   \nonumber \\
						      & \geq & \eta - 2 \epsilon(1,4,0) \nonumber \\
						      &  >   & \epsilon(1,4,0).      \label{eq:NicePath-1}
\end{eqnarray}
\begin{figure}[ht]
\begin{center}
\psfragscanon{
\psfrag{$p$}{$p$}
\psfrag{$q_1$}{$q_1$}
\psfrag{$s_{k-1}$}{$s_{k-1}$}
\psfrag{$s_{k}$}{$s_{k}$}
\psfrag{$r_{k}$}{$r_{k}$}
\psfrag{$u$}{$u$}
\psfrag{$v$}{$v$}
\includegraphics[width=0.6\textwidth]{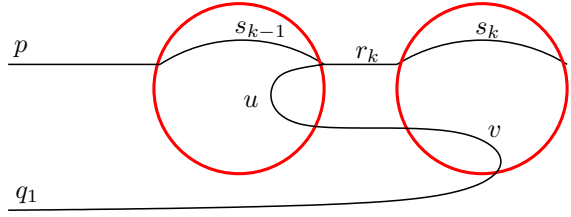}
}
\end{center}
\caption{ Claim 1 in Lemma~\ref{lemma:mainproposition}. If $q_1$ has a $\mathcal{H}$-component $v$ connected to $s_k$, then $dist_X(u_-, u_+) \leq \epsilon(1,4,0)$,
which is a contradiction.}
\label{fig:LemmaNicePath-2}
\end{figure}
Suppose $q_1$ has a $\mathcal{H}$-component $v$ connected to $s_k$. (See Figure~\ref{fig:LemmaNicePath-2}.) 
Consider the subpath $[v_+, u_-]_{q_1}$ of $q_1$ from $v_+$ to $u_-$.
Corollary~\ref{cor:TubeLemma} implies that $[v_+, u_-]_{q_1}$ and $r_k^{-1}$
are $\epsilon(1,4,0)$-similar paths. In particular 
$dist_X(u_-, u_+) \leq \epsilon(1,4,0)$ which contradicts ~\eqref{eq:NicePath-1}. 

{\it Claim 2.} $q$ has a $\mathcal{H}$-component $t_k$ connected to $s_k$. \newline
Consider the triangle $\Delta$ whose sides are $q_1$, $q_2$, and $q$.
By induction hypothesis, $q_2$ has a $\mathcal{H}$-component $t$ which is connected and $\epsilon(1,4,0)$-similar to $s_k$.
By the triangle inequality and our choice of $\eta$, 
\begin{equation}\label{eq:NicePath-2} dist_X(t_-, t_+) \geq \eta - 2 \epsilon(1,4,0) > \tau, \end{equation}
where $\tau$ is the constant of Corollary~\ref{cor:TraingleBigThm}.  Since $q_1$ has no $\mathcal{H}$-component connected to $s_k$, 
Corollary~\ref{cor:TraingleBigThm} implies that $q$ has a $\mathcal{H}$-component $t_k$.

{\it Claim 3.} $q$ has a $\mathcal{H}$-component $t_{k-1}$ connected to $s_{k-1}$. \newline
Since $p_2$ and $q_2$ are $0$-similar, if $q_2$ has a $\mathcal{H}$-component $t$ connected to $s_{k-1}$, 
then the BCP-property shows that $dist_X(t_-, t_+) \leq \epsilon(1,2,0)$.
By the induction hypothesis, the $\mathcal{H}$-component $u$ of $q_1$ connected to $s_{k-1}$ 
satisfies $dist_X(u_-, u_+) \geq \eta - 2 \epsilon(1,4,0)$.
Consider the triangle $\Delta$ whose sides are $q_1$, $q_2$, and $q$.
Corollary~\ref{cor:TraingleBigThm} implies that $q$ has a $\mathcal{H}$-component $t_{k-1}$ connected to $s_{k-1}$ such that
\begin{equation}\label{eq:NicePath-3}  dist_X( (t_{k-1})_- , (t_{k-1})_+) \geq \eta - 2\epsilon(1,4,0) - \epsilon(1,2,0) - 3\tau > 0 .\end{equation}
The last inequality follows by our choice of $\eta$.

{\it Claim 4.} for each $1 \leq i < k-1$, $q$ contains a $\mathcal{H}$-component $t_i$ which is connected to $s_i$. \newline
Let $q'$ be the subpath of $q$ from $q_-$ to $(t_{k-1})_+$.
Notice that $q'$ and $q_1$ are $\epsilon(1,4,0)$-similar (Apply Corollary~\ref{cor:TubeLemma} to $r_k$
and the subpath of $q$ from $(t_{k-1})_+$ to $(t_k)_-$).
By induction hypothesis  $q_1$ has a $\mathcal{H}$-component $u_i$ connected to $s_i$ such that
\begin{equation} \label{eq:NicePath-4} dist_X( (u_i)_- , (u_i)_+) \geq \eta - 2 \epsilon(1,2,0) > \epsilon(1,0, \epsilon(1,4,0)) ,\end{equation}
where the last inequality follows by our choice of $\eta$.
Applying the BCP-property to $q'$ and $q_1$, one sees that $q'$ has a $\mathcal{H}$-component $t_i$ connected to $s_i$ such that
\begin{equation} \label{eq:NicePath-5} dist_X( (t_i)_- , (t_i)_+) \geq \eta - 2 \epsilon(1,2,0) - 2\epsilon(1,0, \epsilon(1,4,0)) > 0, \end{equation}
where the last inequality follows by our choice of $\eta$.
The claim follows.

{\it Claim 5.} The $\mathcal{H}$-components $s_i$ and $t_i$ of $p$ and $q$, respectively, are
$\epsilon(1,4,0)$-similar. \newline
Let $w_i$ be the subpath of $q$ between $t_i$ and $t_{i+1}$.
Corollary~\ref{cor:TubeLemma} implies that $r_i$ and $w_i$ are
$\epsilon(1,4,0)$-similar for each $i$. This completes the proof of the lemma.  
\end{proof}

\begin{lemma}
Let $p$ be a path satisfying the hypothesis of Proposition ~\ref{prop:nicepath}.
Then $p$ is a $(3,0)$-quasi-geodesic. 
\end{lemma}
\begin{proof}
Let $p'$ be a subpath of $p$ and let $q'$ a geodesic in $\GammaH$ connecting the endpoints of $p'$.
The Lemma above implies that $q'$ has a $\mathcal{H}$-component $t_i$ connected to the
$\mathcal{H}$-component $s_i$ of $p'$. Therefore we have decompositions of the two paths of the form:
\begin{eqnarray*} p' = r_\imath '  s_\imath  r_{\imath +1}    \dots  r_{\imath+ \jmath}  s_{\imath+ \jmath}  r_{\imath+ \jmath + 1}' \end{eqnarray*}
\begin{eqnarray*} q' = u_\imath  t_\imath   u_{\imath +1}    \dots  u_{\imath+\jmath}  t_{\imath + \jmath}  u_{\imath+ \jmath + 1}  .\end{eqnarray*}
Corollary~\ref{cor:TubeLemma} implies that $l(r_i) \leq l(u_i) + 2$ for $i=\imath+1 \dots \imath + \jmath$,
$l(r_\imath') \leq l(u_\imath) +2$, and $l(r_{\imath+\jmath+1}')\leq l(u_{\imath+\jmath+1})$. It follows that
\begin{eqnarray*}
 l(p') & \leq & l(u_{\imath} ) +2 +  \sum_{i=\imath+1}^{\jmath+\jmath} \big ( l(s_i) + l(u_{i+1}) +2 \big )  +    l(u_{\imath + \jmath + 1}) +2 \\
        & \leq & 3 l(q') 
\end{eqnarray*}
Since $p'$ was arbitrary, it follows that $p$ is a $(3,0)$-quasi-geodesic in \newline $\GammaH$.
\end{proof}

\begin{lemma}
Let $p$ be a path satisfying the hypothesis of Proposition ~\ref{prop:nicepath}.
Then $p$ is without backtracking. 
\end{lemma}
\begin{proof}
Suppose $p$ backtracks. Let $u$ and $v$ be different connected $\mathcal{H}$-components of $p$ such that the subpath $r$ in-between is without backtracking. 
Observe that $r$ contains one of the $\mathcal{H}$-components $s_i$ of $p$. 
Since $r$ is a $(3,0)$-quasi-geodesic, the BCP-property implies \[ dist_X( (s_i)_-, (s_i)_+ ) \leq \epsilon(3,0,0).\] 
But this is a contradiction since 
\begin{equation} \label{eq:NicePath-6}  dist_X( (s_i)_-, (s_i)_+ ) \geq \eta > \epsilon(3,0,0) , \end{equation}
by our choice of $\eta$.
\end{proof}

A lower bound for $\eta$ is given by  
~\eqref{eq:NicePath-1},   ~\eqref{eq:NicePath-2}, 
~\eqref{eq:NicePath-3},    ~\eqref{eq:NicePath-4},    ~\eqref{eq:NicePath-5}, and  ~\eqref{eq:NicePath-6}.

\section{Quasiconvex Subgroups}\label{sec.quasiconvexity}
This section consist of three parts. 
First (relatively) quasiconvex subgroups are defined following D.\ Osin's work in ~\cite{Os06}. 
The second part consists of the proof of Proposition ~\ref{prop:IntersectionQuasiconvexGroups} which states that  the collection of quasiconvex 
subgroups of a relatively hyperbolic group is closed  under finite intersections, and the third part consists of the proof of 
Proposition ~\ref{prop:ParabolicClasses} on maximal parabolic subgroups of quasiconvex subgroups.

In this section, $G$ is a group generated by a finite set $X$, $\Hs$ a collection of subgroups of $G$, and $G$ is hyperbolic relative to  $\Hs$.

\subsection{Definition of Relatively Quasiconvex Subgroups}
Relatively  quasiconvex subgroups of relatively hyperbolic groups were introduced by D. Osin in ~\cite{Os06} as a generalization of quasiconvex subgroups of word-hyperbolic groups.
F.\ Dahmani in ~\cite{DA03} introduced a dynamical definition of quasiconvex subgroups in relatively hyperbolic groups. C.\ Hruska  showed that both notions are equivalent ~\cite{HK08}..

\begin{definition}\label{defn:Quasi-convexSubgroup}[D.\ Osin ~\cite{Os06}]    A subgroup $Q$ of $G$ is called {\it quasiconvex relative to $\Hs$}
(or simply quasiconvex when the collection $\Hs$ is fixed) if there exists a constant $\sigma \geq 0$ such that the following condition holds. 
Let $f$, $g$ be two elements of $Q$, and $p$ an arbitrary geodesic path from $f$ to $g$ in $\GammaH$. Then for any vertex $v \in p$, there
exists a vertex $w \in Q$ such that $ dist_X(u,w) \leq \sigma .$
\end{definition}

\subsection{Proof of Proposition ~\ref{prop:IntersectionQuasiconvexGroups}}

The following lemma is used several times in the paper, in particular in the  proof of Proposition ~\ref{prop:IntersectionQuasiconvexGroups}. 
In the context of countable groups, we said that a left invariant metric is \emph{proper} if balls of finite radius contain a finite number of elements.
\begin{lemma}\label{lem:quasiorthogonality1}
Let $A$ be a countable group with a proper left invariant metric $d$.
Then for any subgroups $B$ and $C$ of $A$, and any constant $K \geq 0$,
there exists $M = M(B,C,K) \geq 0$ so that
\[  B \cap N_{K}(C)  \subset  N_M ( B \cap C), \]
where $N_{K}(C)$ and $N_M (B \cap C)$ denote the closed $K$-neighborhood and the closed $M$-neighborhood
of $C$ and $B\cap C$ in $(A,d)$ respectively.
\end{lemma}
\begin{proof}
Suppose the statement is false for the constant $K$. Then there are sequences $\{ q_n \}_{n=1}^\infty$
and $\{ h_n \}_{n=1}^\infty$ such that $q_n \in B$, $q_n h_n \in C$, $d(1,h_n) \leq K$, and
\[ d(q_n, B\cap C ) \geq n .\]
Since balls are finite in the metric space $(A,d)$, without lost of generality
assume $\{ h_n \}_{n=1}^\infty$ is a constant sequence $\{ h\}_{n=1}^\infty$. 
For any $m $ and $n$, observe that $q_nq_m^{-1}
= (q_nh) (q_m h)^{-1} \in B\cap C$, and hence $q_mh$ and $q_nh$
are in the same right coset of $B\cap C$, say $ (B\cap C) f$. It
follows that
\[ d(q_n, B \cap C ) \leq d(q_n, q_nh) + d(q_nh, B \cap C) \leq K + d(1,f) \]
for any $n$, a contradiction.
\end{proof}
\begin{remark}
A more general result than Lemma~\ref{lem:quasiorthogonality1} appears in ~\cite[Proposition 9.4]{HK08}.
\end{remark}

\begin{proof}[Proof of Proposition ~\ref{prop:IntersectionQuasiconvexGroups}]
Let $\sigma>0$ so that $Q$ and $R$ are $\sigma$-quasiconvex relative to $\Hs$. Since the generating set $X$ is finite, the metric $dist_X$ on $G$ is proper.
Let $M=M(Q,R,2\sigma)$ be the constant given by Lemma~\ref{lem:quasiorthogonality1}
satisfying
\begin{eqnarray*}  Q \cap N_{2\sigma}(R)  \subset  N_M (Q \cap R) ,  \end{eqnarray*}
where the neighborhoods are taken in the metric space $(G,dist_X)$.

We claim $Q\cap R$ is a $(\sigma+M)$-quasiconvex relative to $\Hs$.
Let $g \in Q\cap R$, let $p$ be a geodesic from $1$ to $g$ in $\GammaH$, and let $u$ be a vertex of $p$.
Since $Q$ and $R$ are $\sigma$-quasiconvex, there exists $s \in Q$ and $t \in R$ so that
\[ max\{ dist_X(s,u), dist_X(t,u) \} \leq \sigma .\]
It follows that $s \in Q \cap N_{2\sigma}(R) $, and hence there is $v \in Q \cap R$ so that
$ dist_X(s, v) \leq M .$ Therefore $v \in Q\cap R$ and  $dist_X(u,v)  \leq \sigma+M.$
\end{proof}

\subsection{Proof of Proposition ~\ref{prop:ParabolicClasses}}

\begin{propositionB}
Let $Q$ be a $\sigma$-quasiconvex subgroup of $G$. Then any infinite maximal parabolic subgroup of $Q$ is conjugate by an element of $Q$ to a subgroup in the set
\begin{equation*} \{  Q \cap H^z  :  H \in \mathcal{H} \text{ and }  z \in G \text{ with } |z|_X \leq \sigma \}. \end{equation*}
In particular, the number of infinite maximal parabolic subgroups up to conjugacy in $Q$ is finite.
\end{propositionB}
\begin{proof}
Let $g \in G$ and $H \in \mathcal{H}$. Suppose that $Q \cap H^g$ is an infinite subgroup. 
Since the generating set $X$ of $G$ is finite, there is an element $h \in H$ such that $|h|_X > \epsilon(1, 0, |g|_X)$
and $h^g \in Q \cap H^g$.  Let $p$ be a geodesic from $1$ to $h^g$. Then the BCP-property~\ref{BCP} implies that $p$ has an 
$\mathcal{H}$-component $s$ contained in the left coset $gH$. Since $Q$ is $\sigma$-quasiconvex, there is an element $y \in Q$ such that $dist_X (y, s_-) \leq \sigma$. The group element $z = y^{-1}s_-$
satisfies that $|z|_X \leq \sigma$ and 
\[ (Q \cap H^z)^y =  Q \cap H^{s_-} = Q\cap H^g . \]
\end{proof}

\section{Proofs of the Combination Theorems and Applications}\label{sec.mainproofs}

The proofs of Theorems~\ref{thm:main_intro} and~\ref{thm:main_intro_2}  adapt some of Gromov's ideas in ~\cite[section 5.3.C.]{Gr87} on 
combination theorems for quasiconvex subgroups in word-hyperbolic groups. We sketch the general argument. Suppose $Q\ast_{Q\cap R} R$ is an amalgamated
product of quasiconvex subgroups of $G$ satisfying the conditions of one of  main theorems. Given a non-trivial element $f$ of
$Q\ast_{Q\cap R} R$, we use its normal form to produce a path $o$ in the relative Cayley graph of $G$ from $1$ to the image of $f$.
Then the path $o$ is shortened by replacing each $\mathcal{H}$-component with more than one edge by a single edge;
the new path is denoted by $p$. (See Figure~\ref{fig:p_thm1_1} below.)  Proposition ~\ref{prop:nicepath} implies that $p$ is a
$(\lambda_0,0)$-quasi-geodesic with different endpoints, proving that the map from $Q\ast_{Q\cap S} S$ into $G$ is injective. Since
$\lambda_0$ is independent of the element $f$, the image of $Q\ast_{Q\cap S} S$ in $G$ will be a quasiconvex subgroup.

The section consists of four parts. In the first subsection a proof of Theorem~\ref{thm:main_intro} is explained in detail. 
Then the proof of Theorem~\ref{thm:main_intro_2} is discussed. The last two subsections correspond to the proof of Corollaries~\ref{cor:double_quasiconvex} and~\ref{cor:fully_quasiconvex}.

\subsection{Proof of Theorem~\ref{thm:main_intro}}
Conjugate the subgroup $Q$ if necessary and assume that $P = H$ for some $H \in \mathcal{H}$.
A lower bound for the constant $C$ is defined during the course of the proof, in particular, the constant $C$ is
chosen large enough to satisfy ~\eqref{eq:maintheorem-1} below. The proof consists of four lemmas. 
Let $\sigma$ be the quasiconvexity constant for $Q$.
\begin{lemma}\label{lem:main_theorem_1}
The natural homomorphism $Q \ast_{Q\cap R} R \longrightarrow G$ is injective.
\end{lemma}
\begin{proof}
Let $f$ be a non-trivial element of $Q \ast_{Q\cap R} R$.
If $f$ is conjugate to an element of $Q$ or an element of $R$,
then it is clear that its image is not trivial.
Otherwise,
\begin{equation}\label{eq:normalform} f = g_1h_1 \dots g_kh_k  \end{equation}
where  $g_i \in Q\setminus Q\cap R$ for $1< i \leq k$,  $h_i \in R \setminus Q\cap R$ for $1\leq i < k$,
either $g_1=1$ or $g_1 \in Q\setminus Q\cap R$, and either $h_k=1$ or $h_k \in R \setminus Q\cap R$.
Since $f$ is not conjugate to an element of $Q$ or $R$, after conjugating if necessary, assume that $g_1\not =1$ and $h_k \not = 1$.
Consider the path $o$ in $\GammaH$ given by
\begin{eqnarray*}\label{eq:path_o} o = u_1  v_1  \dots  u_k  v_k   \end{eqnarray*}
where each $u_i$ and $v_i$ are geodesic paths in $\GammaH$,  
$Label(u_i)$ represents $g_i$, and $Label(v_i)$ represents $h_i$.

{\it Claim 1.} Let $t_i$ be the $\mathcal{H}$-component of the path $o$ containing the
subpath $v_i$. Then the $\mathcal{H}$-components $t_i$ and $t_{i+1}$ are not connected.
\newline If $t_i$ and $t_{i+1}$ are connected,
then $Label(u_i)$ represents an element of $R$. But this contradicts the assumptions on
the normal form ~\eqref{eq:normalform} of the element $f$.

\begin{figure}[ht]
\begin{center}
\psfragscanon{
\psfrag{$x_1$}{$x_1$}
\psfrag{$x_2$}{$x_2$}
\psfrag{$x_3$}{$x_3$}
\psfrag{$x_4$}{$x_4$}
\psfrag{$x_5$}{$x_5$}
\psfrag{$x_6$}{$x_6$}
\psfrag{$u_i$}{$u_i$}
\psfrag{$u_{i+1}$}{$u_{i+1}$}
\psfrag{$z_1$}{$z_1$}
\psfrag{$z_2$}{$z_2$}
\includegraphics[width=0.7\textwidth]{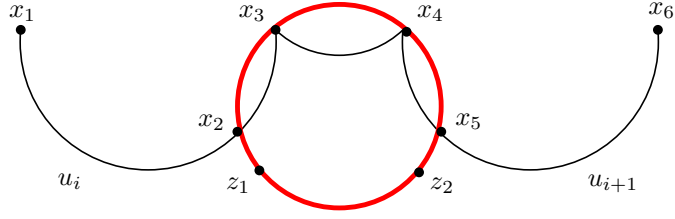}
}
\end{center}
\caption[Proof of Theorem~\ref{thm:main_intro}]{ The endpoints of
$v_i$ are $x_3$ and $x_4$. The $\mathcal{H}$-component $t_i$ of
$o$ containing $v_i$ consists of the path $[x_2,x_3] v_i
 [x_4,x_5]$. } \label{fig:p_thm1_2}
\end{figure}

{\it Claim 2.}  For each $i$, $dist_X( (t_i)_-, (t_i)_+ ) > \eta$, where $\eta$ is the constant from Proposition ~\ref{prop:nicepath}.
\newline Fix $i$ and let $x_1=(u_i)_-$, $x_2=(t_i)_-$, $x_3=(u_i)_+$, $x_4=(u_{i+1})_-$, $x_5=(t_i)+$, and $x_6=(u_{i+1})_+$.
(See Figure~\ref{fig:p_thm1_2}.)

First we show that there are elements $z_1$ and $z_2$ in the left coset $x_3 (H\cap Q) = x_4(H\cap Q)$ such that
\begin{eqnarray*} max\{ dist_X(x_2 , z_1) , dist_X(x_5, z_2)  \} \leq  M(H,Q,\sigma) ,  \end{eqnarray*}
where $ M(H,Q,\sigma)$ is the constant provided by Lemma~\ref{lem:quasiorthogonality1} for the subgroups $H$ and $Q$, the
constant $\sigma $, and the proper metric $dist_X$. Since $Q$ is $\sigma$-quasiconvex and $Label(u_i)$ represents an element of
$Q$, we have that $x_3^{-1}x_2 \in H$ and $d_X(x_3^{-1}x_2, Q) \leq \sigma$; then Lemma~\ref{lem:quasiorthogonality1} implies
that
\begin{eqnarray*} d_X( x_3^{-1}x_2, H\cap Q) \leq  M(H,Q,\sigma) ; \end{eqnarray*}
hence there is an element $z_1 \in x_3(H\cap Q)$ such that $d_X(x_2,z_1)\leq M(H,Q,\sigma)$.
A similar argument guarantees the existence of an element $z_2$ with the desire properties.

Since $x_3^{-1}x_4 \in R \setminus Q$, $x_3^{-1}z_1 \in (H\cap Q) \subset R$, and $x_4^{-1}z_2 \in (H\cap Q) \subset R$, we have that
\begin{eqnarray*} z_1^{-1}z_2 = (z_1^{-1}x_3)(x_3^{-1}x_4)(x_4^{-1}z_2) \in R \setminus Q,  \end{eqnarray*}
and hence $dist_X(z_1, z_2) \geq C $, by hypothesis (2) on the length of the elements of $R \setminus Q$. Finally, by our choice
of the constant $C$,
\begin{eqnarray}
  dist_X(x_2, x_5) & \geq & dist_X(z_1, z_2) - dist_X(x_2, z_1) - dist_X(x_5, z_2) \nonumber \\
												  & \geq & C - 2 M(H, Q, \sigma ) \nonumber \\
\label{eq:maintheorem-1}                & \geq & \eta
\end{eqnarray}
which proves the claim.

\begin{figure}[ht]
\begin{center}
\psfragscanon{
\psfrag{$o$}{$o$}
\psfrag{$u_{i-1}$}{$u_{i-1}$}
\psfrag{$v_{i-1}$}{$v_{i-1}$}
\psfrag{$u_i$}{$u_i$}
\psfrag{$v_i$}{$v_i$}
\psfrag{$u_{i+1}$}{$u_{i+1}$}
\psfrag{$v_{i+1}$}{$v_{i+1}$}
\psfrag{$g_iH'$}{$g_iH'$}
\psfrag{$g_{i+1}H''$}{$g_{i+1}H''$}
\includegraphics[width=0.95\textwidth ]{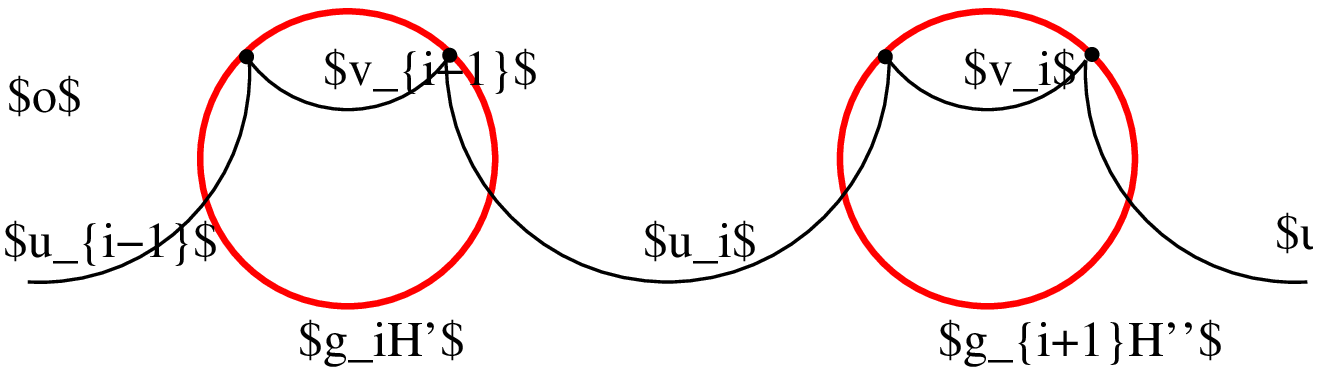}
}
\end{center}

\begin{center}
\psfragscanon{
\psfrag{$p$}{$p$}
\psfrag{$r_{i-1}$}{$r_{i-1}$}
\psfrag{$s_{i-1}$}{$s_{i-1}$}
\psfrag{$r_i$}{$r_i$}
\psfrag{$s_{i}$}{$s_{i}$}
\psfrag{$r_{i+1}$}{$r_{i+1}$}
\psfrag{$s_{i+1}$}{$s_{i+1}$}
\psfrag{$g_i H'$}{$g_i H'$}
\psfrag{$g_{i+1} H''$}{$g_{i+1} H''$}
\includegraphics[width=0.95\textwidth ]{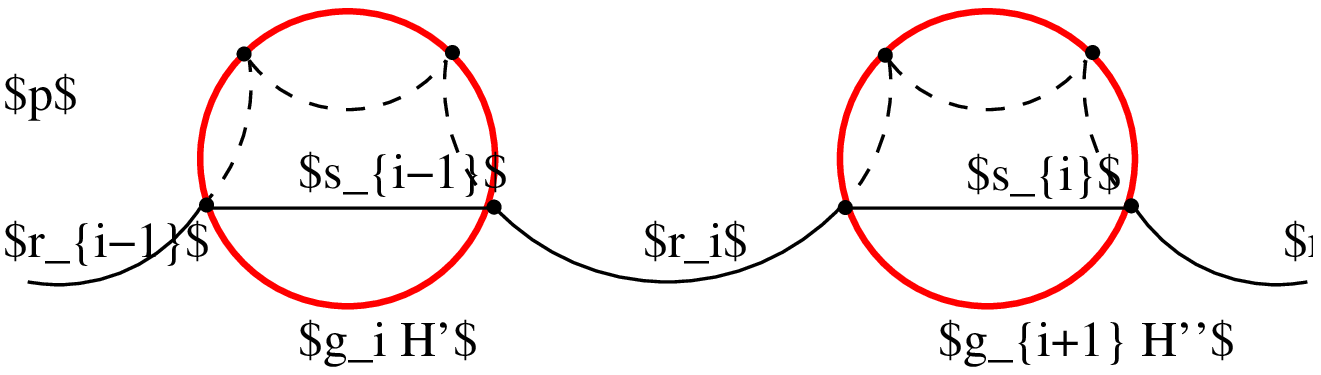}
}
\end{center}
\caption[Proof of Theorem~\ref{thm:main_intro}]{ The path $o$ and the resulting quasi-geodesic $p$.}\label{fig:p_thm1_1}
\end{figure}

{\it Claim 3.} Let $p$ be the path obtained by replacing each $\mathcal{H}$-component $t_i$ of $o$
for a single edge $s_i$. (See Figure~\ref{fig:p_thm1_1}.) Then $p$ is $(\lambda_0,0)$-quasi-geodesic
with different endpoints. In particular the image of $f$ by the map $Q \ast_{Q\cap R} R \longrightarrow G$ is not trivial.
\newline  The path $p$ has a decomposition of the form
\begin{equation} p = r_1  s_1  \dots  r_k  s_k , \end{equation}
where $r_1$ or $s_k$ may be trivial.
The definition of the path $p$ (and the path $o$) show that $r_i$ and $s_i$ are geodesic segments
in $\GammaH$; claim 1 shows that the $\mathcal{H}$-components $s_i$ and $s_{i+1}$ of $p$ are not connected;
and claim 2 implies
\begin{eqnarray*} dist_X( (s_i)_- , (s_i)_+) \geq \eta .  \end{eqnarray*}
Proposition ~\ref{prop:nicepath} implies that $p$ is a $(\lambda_0,0)$-quasi-geodesic
with different endpoints, proving the claim.
\end{proof}
 
\begin{lemma}\label{lemma:main3}
The subgroup $\langle Q \cup R \rangle$ is relatively quasiconvex and its quasiconvexity constant is independent of the choice of $R$.
\end{lemma}
\begin{proof}
Let $f \in \langle Q \cup R \rangle$, and  let $q$ be a geodesic in $\GammaH$ from $1$ to $f$.
If $f \in Q \cup R$ then it is trivial that any vertex of $q$ is at most $\sigma$ apart from an element of $\langle Q \cup R \rangle$ with respect to the metric $dist_X$.
Otherwise, let $p$ be the $(\lambda_0,0)$-quasi-geodesic constructed during the proof of Lemma~\ref{lem:main_theorem_1} from $1$ to $f$.
Notice that any vertex of $p$ is  at most  $\sigma$ apart from an element of $\langle Q \cup R \rangle$ with respect to the metric $dist_X$.
The BCP property implies that every vertex of $q$ is at most $\epsilon(\lambda_0,0,0)$ apart from the set of vertices of $p$. 
It follows that any vertex of $q$ is at most $(\sigma + \epsilon(\lambda_0,0,0))$ apart from $\langle Q \cup R \rangle$.
This shows $\langle Q \cup R \rangle$ is relatively $(\sigma + \epsilon(\lambda_0,0,0))$-quasiconvex.
\end{proof}
\begin{lemma} \label{lemma:main2} Any parabolic element of the subgroup $\langle Q\cup R \rangle$ is either
conjugate to an element of $Q$ or to an element of $R$ by an element of $\langle Q\cup R \rangle$.
\end{lemma}
\begin{proof}
If $f$ is a parabolic element of $G$, then its action on $\GammaH$ fixes setwise a subset of diameter one. 
Indeed, if $f$ is a parabolic element, then there is an element $g \in G$ such that $gfg^{-1} \in H_i$ for some $i$. 
It follows that $f$ fixes setwise the left coset $g^{-1}H_i$ which has diameter one.

Let $f$ be an element of $Q \ast_{Q\cap R} R$, and suppose that $f$ is not conjugate to an element of $Q$ or $R$. 
We claim that $f$ acts on a bi-infinite quasi-geodesic $\widetilde{p}$ in $\GammaH$ as a non-trivial translation, and hence the observation of the previous paragraph implies that $f$ is not a parabolic element. 

Conjugate $f$, if necessary, and assume that its normal form,
\begin{eqnarray*} f = g_1h_1 \dots g_kh_k,  \end{eqnarray*}
satisfies $g_1 \not = 1$ and $h_k \not = 1$.  Consider the path $o$ in $\GammaH$
\begin{eqnarray*} o \  = \ u_1  v_1  \dots  u_k  v_k \end{eqnarray*}
from $1$ to $f$, where each $u_i$ and $v_i$ are non-trivial geodesics in $\GammaH$,  $Label(u_i)$ represents $g_i$,
and $Label(v_i)$ represents $h_i$.  Let $\widetilde{o}$ the bi-infinite path
\begin{eqnarray*} \widetilde{o} \  =   \   \dots  f^{-3}(o)\   f^{-2}(o) \   f^{-1}(o)\   o\   f(o)\   f^2(o)\   f^3(o)   \dots ,  \end{eqnarray*}
and let $\widetilde{p}$ the path obtained by replacing each $\mathcal{H}$-component of $\widetilde{o}$ by a
single edge. The argument of Lemma~\ref{lem:main_theorem_1} shows that the subpath of $\widetilde{p}$
induced by the subpath 
\begin{eqnarray*}    f^{-k}(o)\   f^{-k+1}(o)  \dots  f^{-1}(o)\  o\  f(o)  \dots  f^{k-1}(o) \  f^k(o),  \end{eqnarray*}
of $\widetilde{o}$ is a $(\lambda_0,0)$-quasi-geodesic for any integer $k>0$. It follows that $\widetilde{p}$ is a bi-infinite $(\lambda_0,0)$-quasi-geodesic, and that
the (image in $G$ of the) element $f$ acts as a nontrivial translation on this bi-infinite quasi-geodesic.
\end{proof}
\begin{lemma}\label{lemma:main4}
Any parabolic subgroup of $\langle Q\cup R \rangle$ is conjugate either to a subgroup of $Q$ or to a subgroup of $R$ by an element of $\langle Q\cup R \rangle$.

In particular, if $\{ K_1, \dots, K_l \}$ is the collection of maximal parabolic subgroups of $Q$ up to conjugacy in $Q$. Then
the collection of maximal parabolic subgroups of $\langle Q\cup R\rangle$ up to conjugacy in $\langle Q\cup R\rangle$ is
\begin{enumerate}
\item $\{R, K_1, \dots , K_l\}$, if $Q \cap P$ is trivial;
\item $\{R, K_2, \dots , K_l\}$, if $Q\cap P$ and $K_1$ are conjugate in $Q$.
\end{enumerate}
\end{lemma}
\begin{proof}
Assume that $Q\cap R$ is a proper subgroup of $R$ and $Q$; otherwise there is nothing to prove. 

An easy argument using normal forms shows that if $J$ is a subgroup of $Q\ast_{Q\cap R}R$ that can not be conjugated into $Q$ or $R$, then $J$
contains an element that can not be conjugated into $Q$ or $R$.  By Lemma~\ref{lemma:main2}, any parabolic subgroup $J$ of $Q\ast_{Q\cap R}R$ is conjugate
to a parabolic subgroup of $Q$ or $R$. This also implies that any maximal parabolic subgroup $J$ of $Q\ast_{Q\cap R}R$ is conjugate to a maximal parabolic subgroup of $Q$ or $R$, 
and hence to a subgroup in $\{R, K_1,  K_2, \dots , K_l\}$.

The second statement follows from the following observation. 
For any element $g \in \langle Q\cup R \rangle$, if $g \not \in Q$ then $Q\cap Q^g$ is either trivial or contained in $Q\cap R$;
therefore $K_i$ and $K_j$ are conjugate in $Q\ast_{Q\cap R}R$ only if $i = j$.  Since $Q\cap R$ is a proper subgroup of $R$ and $Q$, $R$ is
not conjugate to a subgroup of $Q$ in $Q\ast_{Q\cap R}R$.
\end{proof}

\subsection{Proof of Theorem~\ref{thm:main_intro_2} }

Conjugate the subgroups $Q_1$ and $Q_2$ if necessary and assume that $P = H$ for some $H \in \mathcal{H}$.
A lower bound for the constant $C$ is given by ~\eqref{eq:maintheorem-2-1} below. 
The proof is completely analogous to the proof of Theorem~\ref{thm:main_intro}. 
Let $\sigma$ be a common quasiconvexity constant for $Q_1$ and $Q_2$.

\begin{lemma}\label{lem:main2-1}
The natural homomorphism $Q_1 \ast_R hQ_2h^{-1} \longrightarrow G $ is injective.
\end{lemma}
\begin{proof}
Let $f$ be an element of $Q_1 \ast_{R} hQ_2h^{-1}$ with normal form 
\begin{equation}\label{eq:normalform2} f = g_1 g_2^h \dots g_{2k-1}g_{2k}^h \end{equation}
where  $g_{2i+1} \in Q_1 \setminus R$ for $1 \leq i \leq k$,
$g_{2i} \in Q_2 \setminus R$ for $1 \leq i < k$,
either $g_1=1$ or $g_1 \in Q_1 \setminus R$, and either $g_{2k}=1$ or $g_{2k} \in Q_2 \setminus R$.
Consider the path $o$ in $\GammaH$ from $1$ to $f$ given by
\begin{eqnarray*} o = u_1   v_1   u_2  v_2    \dots    u_{2k-1}  v_{2k-1}  u_{2k}   v_{2k}  \end{eqnarray*}
where  $u_i$, $v_{2i}$ and $v_{2i+1}$ are geodesic paths  in $\GammaH$ such that $Label(u_i)$ represents $g_i$,
$Label(v_{2i-1})$ represents $h$, and $Label(v_{2i})$ represents $h^{-1}$.

{\it Claim 1.} Let $t_i$ be the $\mathcal{H}$-component of $o$ that contains the subpath $v_i$.
Then the $\mathcal{H}$-components $t_i$ and $t_{i+1}$ are not connected.\newline
If $t_{2i-1}$ and $t_{2i}$ are connected $\mathcal{H}$-components, then
$g_{2i} \in R$, which contradicts the assumptions on
~\eqref{eq:normalform2}. Analogously $t_{2i}$ and $t_{2i+1}$ are not connected $\mathcal{H}$-components.

{\it Claim 2.} $dist_X( (t_i)_-, (t_i)_+ ) > \eta$ for each $i$, where $\eta$ is the constant from Proposition ~\ref{prop:nicepath}. \newline
Fix an odd value of $i$, and let $x_1=(u_i)_-$, $x_2=(t_i)_-$, $x_3=(u_i)_+$, $x_4=(u_{i+1})_-$, $x_5=t_+$
and $x_6=(u_{i+1})_+$. The argument used to prove Claim 2 of Lemma~\ref{lem:main_theorem_1} shows that
there are elements $z_1$ in the left coset $x_3R$, and $z_2$ in the left coset $x_4 R$
such that
\begin{eqnarray*} 			 dist_X(x_2 , z_1)  & \leq  &  M(H,Q_1,\sigma),  \end{eqnarray*}
and
\begin{eqnarray*}			  dist_X(x_5, z_2 )  & \leq  &  M(H,Q_2,\sigma),	  \end{eqnarray*}
where $ M(H, Q_i,\sigma)$ is the constant provided by
Lemma~\ref{lem:quasiorthogonality1} for the subgroups $H$ and $Q_i$, the constant $\sigma $, and the proper metric $d_X$. It follows that
\begin{eqnarray*} z_1^{-1}z_2 = (z_1^{-1}x_3)(x_3^{-1}x_4)(x_4^{-1}z_2) \in RhR ,  \end{eqnarray*}
and hence $ dist_X(z_1, z_2) \geq C .$ Now the triangle inequality and our choice of $C$ implies
\begin{eqnarray}
  dist_X(x_2, x_5) & \geq & dist_X(z_1, z_2) - dist_X(x_2, z_1) - dist_X(x_5, z_2) \nonumber \\
						& \geq & C - M(H, Q_1,\sigma ) - M(H, Q_2,\sigma ) 			    \nonumber \\
\label{eq:maintheorem-2-1}                & \geq  &  \eta,
\end{eqnarray} 
where $\eta$ is the constant from Proposition ~\ref{prop:nicepath}
The case for an even value of $i$ is similar, the only difference is that $z_1^{-1}z_2 \in Rh^{-1}R$, which also implies $ dist_X(z_1, z_2) \geq C$.

{\it Claim 3.} Let $p$ be the path obtained by replacing each $\mathcal{H}$-component $v_i$ of $o$ by a single edge $s_i$. 
The above claims and Proposition ~\ref{prop:nicepath} imply that
 $p$ is $(\lambda_0,0)$-quasi-geodesic with different endpoints. In particular the image of $f$ by the map $Q_1 \ast_R hQ_2h^{-1} \longrightarrow G$
 is not trivial.
 \end{proof}
 
 \begin{remark}[A technical remark on the proof of Lemma~\ref{lem:main2-1}]  \label{rem:additional}
If $p$ is the path from $1$ to an element $f$ of $\langle Q_1, hQ_2h^{-1}\rangle \setminus Q_1 \cup Q_2$ constructed in the proof, 
then $p$ has at least two different $\mathcal{H}$-components $s_1$ and $s_2$ of $X$-length at least $\eta$, namely the ones induced by an element of $hQ_2h^{-1}$ in the normal form of $f$.
Since $\eta$ is larger than $\epsilon (\lambda_0, 0,0)$ (see ~\eqref{eq:NicePath-6}), the BCP-property implies that any geodesic from $1$ to $f$ has 
at least two different $\mathcal{H}$-components.  In particular, the element $f$ does not belong to a subgroup $H \in \mathcal{H}$. 

Therefore, if $Q_1\cap H = Q_2\cap H$ is the trivial subgroup, then $\langle Q_1\cup hQ_2h^{-1} \rangle \cap H$ is trivial.
\end{remark}

Similar arguments to the proofs of Lemmas~\ref{lemma:main2} and~\ref{lemma:main3} show the following.

\begin{lemma}
The subgroup $\langle Q_1\cup hQ_2h^{-1} \rangle$ is relatively quasiconvex.
\end{lemma}

\begin{lemma}\label{lem:main2-3}
Every parabolic subgroup of $\langle Q_1\cup hQ_2h^{-1} \rangle$ is conjugate to a parabolic subgroup of $Q_1$ or $Q_2$ by  an element of $\langle Q_1\cup hQ_2h^{-1} \rangle$.

In particular, if $\mathcal{K}_i$ is the collection of maximal parabolic subgroups of $Q_i$ up to conjugacy in $Q_i$, for $i=1,2$. Then every maximal parabolic
subgroup of $\langle Q_1\cup hQ_2h^{-1} \rangle$ is conjuate to a subgroup in $\mathcal{K}_1 \cup \mathcal{K}_2$ by an element 
of $\langle Q_1\cup hQ_2h^{-1} \rangle$.
\end{lemma}

\subsection{Proof of Corollary~\ref{cor:double_quasiconvex} }

\begin{lemma} \label{lemma:abeliangroups}
Suppose  $A$ is an abelian group with a finite generating set $Y$, $B$ is a subgroup of $A$, and $h \in A$ such that 
\[ rank_{\integers}(B) < rank_{\integers}( \langle B\cup \{h\}\rangle ) .\] 
Then there is a constant $\lambda=\lambda (B,h,Y)$ such that  $|g|_Y \geq  \lambda |j|$ for any $j \in \integers$ and $g \in h^j B$
\end{lemma}
\begin{proof}
Since $h$ has infinite order, one can regard $A$ as the direct product $A_1\oplus \langle h_1 \rangle$ where $h \in \langle h_1 \rangle$ and $B\subset A_1$. 
Suppose $Y$ contains $h_1$, and $Y \setminus \{h_1\}$ generates $A_1$. Then $|g|_Y \geq |h^j|_Y \geq |j|$ for any $g \in h^j B $. Since the word metrics associated to
different finite generating sets are Lipschitz equivalent the result follows.
\end{proof}

Corollary~\ref{cor:double_quasiconvex} follows from the following proposition.

\begin{proposition}
Let $G$ be hyperbolic relative to a collection of free abelian subgroups $\mathcal{H}$, $Q$ a relatively quasiconvex subgroup, and $P$ 
a maximal parabolic subgroup of $G$.  If 
\[rank_{\integers} (Q\cap P) < rank_{\integers} (P),\]
then there is an element $h \in P$ with the following property.
For any positive integer $k$, there exist integers $n_1, \dots , n_k$ such that 
\begin{eqnarray*} \Big\langle \bigcup_{i=1}^k h^{n_i}Qh^{-n_i} \Big\rangle  \cong  \Delta_k (Q, Q\cap H) . \end{eqnarray*}
Moreover, the above subgroup is relatively quasiconvex.
\end{proposition}
\begin{proof}
Let $G$ be hyperbolic relative to a collection of abelian subgroups $\mathcal{H}$, $Q$ a relatively quasiconvex subgroup, and $H \in \mathcal{H}$ such that 
\[rank_{\integers} (Q\cap H) < rank_{\integers} (H).\]

Let $h \in H$ be any element such that
\begin{equation} rank_{\integers} ( Q\cap H) <  rank_{\integers} \Big(  \langle (Q\cap H) \cup \{h \} \rangle \Big) , \end{equation}
and let $Y$ be a finite generating set of $H$. By Lemma~\ref{lemma:abeliangroups} there is a constant $\lambda > 0$  such that
\begin{equation} \label{eq:corollary2}  |g|_Y \geq \lambda |j|, \end{equation}
for any integer $j$, and any element $g$ in the left coset $h^j (Q\cap H)$.

By induction on $k$, we prove the existence of integers $n_1, \dots , n_k$ such that the subgroup
$ \big\langle \bigcup_{i=1}^k h^{n_i}Qh^{-n_i} \big\rangle  $
is quasiconvex, isomorphic to $\Delta_k (Q, Q\cap H) $, and 
$Q\cap H =  \big\langle \bigcup_{i=1}^k h^{n_i}Qh^{-n_i} \big\rangle \cap H.$

The case $k=1$ is trivial taking $n_1=0$. Suppose
$R_{k-1} = \big\langle  \bigcup_{i=1}^{k-1} f^{n_i}Qf^{-n_i} \big\rangle$ is a quasiconvex subgroup
isomorphic to $\Delta_{k-1} (Q, Q\cap H)$, and $R_{k-1}\cap H = Q\cap H$.  
Let $C=C(R_{k-1},Q, H)$ be the constant provided by Theorem~\ref{thm:main_intro_2}, and let $n_{k}$ be any integer such that
\begin{equation} \label{eq:cor2-1}  max\{ |g|_{Y} : g \in H \  \text{with}\  |g|_X < C \} <  \lambda n_{k}.   \end{equation}
Now we show that  the quasiconvex subgroups $R_{k-1}$ and $Q$, the maximal parabolic subgroup $H$, and the element $h^{n_{k+1}} \in H$ satisfy
the hypothesis of Theorem~\ref{thm:main_intro_2}:
first, since $H$ is abelian,  $h^{n_{k}} (Q\cap H) h^{-n_{k}} = Q\cap H$;
second, if $g$ is an element of the left coset $h^{n_{k+1}} (Q\cap H)$, then
~\eqref{eq:corollary2} and ~\eqref{eq:cor2-1} imply  that $ |g|_X \geq  C . $

Therefore Theorem~\ref{thm:main_intro_2} implies 
that the subgroup $R_k=\big\langle R_{k-1}  ,  h^{ n_{k}} Q h^{-n_{k}} \big\rangle $ is isomorphic to
$R_{k-1} \ast_{Q \cap H} Q  \cong \Delta_{k} (Q, Q\cap H)$ and is quasiconvex.
We claim that $Q\cap H= R_k \cap H$.

If $Q\cap H$ is not trivial, then Lemma~\ref{lem:main2-3} applied to $Q$ and $R_{k-1}$ implies that $Q\cap H$ is a maximal parabolic subgroup of $R_k$.
Since $R_k \cap H$ is a maximal parabolic subgroup of $R_k$ containing $Q\cap H$, if follows that $Q\cap H = R_k\cap H$.
If $Q\cap H$ is trivial, then remark~\ref{rem:additional} implies $R_k \cap H$ is trivial.
\end{proof}

\subsection{Proof of Corollary~\ref{cor:fully_quasiconvex}}

\begin{figure}[ht]
\begin{center}
\psfragscanon{
\psfrag{$A_1$}{$A_1$}
\psfrag{$A_2$}{$A_2$}
\psfrag{$A_n$}{$A_n$}
\psfrag{$K_1$}{$K_1$}
\psfrag{$K_2$}{$K_2$}
\psfrag{$K_n $}{$K_n$}
\psfrag{$Q$}{$Q$}
\includegraphics[width=0.32\textwidth]{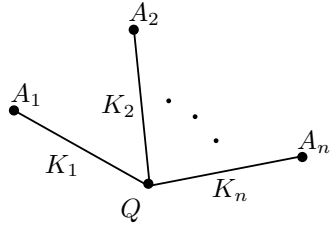}
}
\end{center}
\caption{Tree of groups decomposition of the fully quasiconvex subgroup $R$ of Corollary~\ref{cor:fully_quasiconvex}.}\label{fig:corollaryfullyquasiconvex}
\end{figure}

Corollary~\ref{cor:fully_quasiconvex} follows from the following proposition.

\begin{proposition}
Let $G$ be hyperbolic relative to a collection of free abelian subgroups $\mathcal{H}$, and let $Q$ be a relatively quasiconvex subgroup. 
Then there exists a fully quasiconvex subgroup $R$ which has the tree of groups decomposition described in Figure~\ref{fig:corollaryfullyquasiconvex}, where 
$\{K_1, \dots, K_n\}$ is a collection of representatives of the maximal parabolic subgroups of $Q$ and each $A_i$ is a finite index subgroup of a maximal parabolic subgroup of $G$.
\end{proposition}
\begin{proof}
By Proposition ~\ref{prop:ParabolicClasses},  a collection of representatives of the infinite maximal parabolic subgroups of $Q$ up to conjugacy in $Q$
is finite, say  $K_1, \dots, K_n$.   The desired group is obtained after an $n$-step process which produces a sequence $\{Q_j\}_{j=0}^n$ of quasiconvex subgroups of $G$ where $Q_0=Q$ and $Q_n=R$ is a fully quasiconvex subgroup.
For $1\leq j \leq n$, the group $Q_j$ has the graph of groups decomposition described in Figure~\ref{fig:graph_proof},
where $\{A_1, \dots , A_j, K_{j+1}, \dots , K_n \}$ is the collection of all maximal parabolic
subgroups of $Q_j$ up to conjugation in $Q_j$, and  $A_i$ is a finite index subgroup of a
maximal parabolic subgroup of $G$ for each $i\leq j$. 

Now we explain how to obtain $Q_{i+1}$ from $Q_i$ when $i<n$.
Let $P$ be the maximal parabolic subgroup of $G$ containing $K_{i+1}$, and let $Y$ be a finite generating set of $P$.
Let $C = C(Q_i,P) \geq 0$ the constant provided by Theorem~\ref{thm:main_intro}, and define
\begin{eqnarray*} D = max\{ |g|_{Y} : g \in P \  \text{with}\  |g|_X < C \} .  \end{eqnarray*}
Since $P$ is a finitely generated abelian group, there is a finite index subgroup $A_{i+1}$ of $P$ containing $K_{i+1}$ such that
$|g|_Y \geq D$ for any $g \in A \setminus K$. In particular $|g|_X
\geq C$ for any $g \in A_{i+1} \setminus K_{i+1}$, and hence
Theorem~\ref{thm:main_intro} and Lemma~\ref{lemma:main4}  imply that the subgroup $Q_{i+1} =
\langle Q_i \cup A_{i+1} \rangle $ of $G$ is isomorphic to $Q_i
\ast_{K_{i+1}} A_{i+1}$, is relatively quasiconvex, and $\{A_1,
\dots , A_{j+1}, K_{j+2}, \dots , K_n \}$ is the collection of all
maximal parabolic subgroups of $Q_{i+1}$ up to conjugation in $Q_{i+1}$.  

\begin{figure}[ht]
\begin{center}
\psfragscanon{
\psfrag{$A_1$}{$A_1$}
\psfrag{$A_j$}{$A_j$}
\psfrag{$K_1$}{$K_1$}
\psfrag{$K_j$}{$K_j$}
\psfrag{$K_{j+1}$}{$K_{j+1}$}
\psfrag{$K_n$}{$K_n$}
\psfrag{$Q$}{$Q$}
\includegraphics[width=0.32\textwidth]{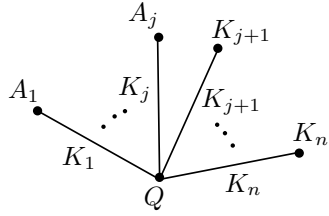}
}
\end{center}
\caption{The graph of groups decomposition of the subgroup $Q_j$ in the proof of Corollary~\ref{cor:fully_quasiconvex}.} \label{fig:graph_proof}
\end{figure}

\end{proof}

\section{Future Directions}\label{sec.questions}
\subsection{Amalgamation along hyperbolic subgroups}
The main results of this paper address Problem~\ref{main.problem}, stated in the introduction, in
the case that $K$ is a maximal parabolic subgroup of $Q_1$ or $Q_2$. The case that $K$ consists
only of hyperbolic elements remains to be considered. 

An analogous version of Problem~\ref{main.problem} for HNN-extensions is of interest.
\begin{problem}
If $R_1$ and $R_2$ are subgroups of a relatively quasiconvex subgroup $Q$ of $G$, 
investigate conditions guaranteeing the existence of a group homomorphism 
$\varphi : R_1 \longrightarrow R_2 $ and an injective homomorphism
\[
Q\ast_{\varphi} \longrightarrow G
\]
with image a quasiconvex subgroup. Here $Q\ast_{\varphi}$ represents the HNN-extension
\[
\presentation{Q, t}{ trt^{-1} = \varphi(r), \ \ r \in R_1 }
\]
\end{problem}

\subsection{Surface subgroups in (relatively) hyperbolic groups.}

Finding subgroups isomorphic to fundamental groups of hyperbolic closed surfaces  in  (relatively) hyperbolic groups has been an theme
in Geometric Group Theory ~\cite{B04, S07}.  D.Cooper, D.Long and A. Reid have produced closed surface subgroups
in fundamental groups of complete hyperbolic 3-manifolds with cusps  ~\cite{CL01, CLR97}.
In ~\cite{BFM07}, N. Brady, M. Forester and the author has shown that a class of word-hyperbolic groups
of the form $F_k \ast_{\integers} F_l$, where $F_k$ and $F_l$ are free groups of rank $k$ and $l$
respectively, have surfaces subgroups.
The techniques used in ~\cite{BFM07} resemble Cooper-Long-Reid ideas of doubling surfaces with one boundary component through a combination theorem.
\begin{problem}
For which other classes of word-hyperbolic groups can these ideas produce surface subgroups?
\end{problem}
D.Cooper, D.Long have produced surface subgroups in Dehn fillings of hyperbolic manifolds~\cite{CL01}. 
In the context of hyperbolic Dehn fillings on relatively hyperbolic groups ~\cite{GM06, Os06-1}:
\begin{problem}
Explore the existence of surface subgroups in one-ended word-hyperbolic groups arising as hyperbolic Dehn fillings of relatively hyperbolic groups.
\end{problem}

\bibliographystyle{amsplain} 
\bibliography{xbib}

\end{document}